\documentclass[12pt,a4paper]{article}
\usepackage{amsfonts}
\usepackage{amsthm}
\usepackage{amssymb}
\usepackage[centertags]{amsmath}
\usepackage{amssymb}
\usepackage{enumerate}
\usepackage{graphicx}
\providecommand{\U}[2]{\protect\rule{1.5in}{1.5in}}
\theoremstyle{plain}
\newtheorem{theorem}{Theorem}[section]

\newtheorem{definition}[theorem]{Definition}
\newtheorem{example}[theorem]{Example}
\newtheorem{lemma}[theorem]{Lemma}

\newtheorem{remarks}[theorem]{Remark}

\numberwithin{equation}{section}
\begin{document}
\title{Solving fuzzy two-point boundary value problem using fuzzy Laplace transform} 
\author{Latif Ahmad\footnote{1. Shaheed Benazir Bhutto University, Sheringal, Upper Dir, Khyber Pakhtunkhwa, Pakistan. 2. Institute of Pure \& Applied Mathematics, University of Peshawar, 25120, Khyber Pakhtunkhwa, Pakistan. E-mail: ahmad49960@yahoo.com}, \mbox{} Muhammad Farooq\footnote{ Institute of Pure \& Applied Mathematics, University of Peshawar, 25120, Khyber Pakhtunkhwa, Pakistan. E-mail: mfarooq@upesh.edu.pk}, \mbox{} Saif Ullah\footnote{Institute of Pure \& Applied Mathematics, University of Peshawar, 25120, Khyber Pakhtunkhwa, Pakistan. E-mail: saifullah.maths@gmail.com}, \mbox{} Saleem Abdullah\footnote{Department of Mathematics, Quaid-i-Azam University, Islamabad, Pakistan. E-mail: saleemabdullah81@yahoo.com}}
\bibliographystyle{ams}
\maketitle
\begin{abstract}
A natural way to model dynamic systems under uncertainty is to use fuzzy boundary value problems (FBVPs) and related uncertain systems. In this paper we use fuzzy Laplace transform to find the solution of two-point boundary value under generalized Hukuhara differentiability. We illustrate the method for the solution of the well known two-point boundary value problem $Schr\ddot{o}dinger$ equation, and homogeneous boundary value problem. Consequently, we investigate the solutions of FBVPs under as a new application of fuzzy Laplace transform.
\end{abstract}
\section{Introduction}
``The theory of fuzzy differential equations (FDEs) has attracted much attention in recent years because this theory represents a natural way to model dynamical systems under uncertainty'',  Jamshidi and Avazpour \cite{1a}.
 The concept of fuzzy set was introduce by Zadeh in $1965$ \cite{1}. The derivative of fuzzy-valued function was introduced by Chang and Zadeh in $1972$  \cite{2}. The integration of fuzzy valued function is presented in \cite{8}. Kaleva and Seikala presented fuzzy differential equations (FDEs) in \cite{3,4}. Many authors discussed the applications of FDEs in \cite{5,6,7}. Two-point boundary value problem is investigated in \cite{9}. In case of Hukuhara derivative the funding Green's function helps to find the solution of boundary value problem of first order linear fuzzy differential equations with impulses \cite{10}. Wintner-type and superlinear-type results for fuzzy initial value problems (FIVPs) and fuzzy boundary value problems (FBVPs) are presented in \cite{11}. The solution of FBVPs must be a fuzzy-valued function under the Hukuhara derivative \cite{12,13,14,15,16}. Also two-point boundary value problem (BVP) is equivalent to fuzzy integral equation \cite{17}. Recently in \cite{18,19,20} the fuzzy Laplace transform is applied to find the analytical solution of FIVPs. According to \cite{21} the fuzzy solution is different from the crisp solution as presented in \cite{12,13,14,15,22,23}. In \cite{21} they solved the $Schr\ddot{o}dinger$ equation with fuzzy boundary conditions. Further in \cite{18} it was discussed that under what conditions the fuzzy Laplace transform (FLT) can be applied to FIVPs. For two-point BVP some of the analytical methods are illustrated in \cite{21,24,25} while some of the numerical methods are presented in \cite{1a,27}. But every method has its own advantages and disadvantages for the solution of such types of fuzzy differential equation (FDE). In this paper we are going to apply the FLT on two-point BVP \cite{21}. Moreover we investigate the solution of second order $Schr\ddot{o}dinger$ equation and other homogeneous boundary value problems \cite{21}. After applying the FLT to BVP we replace one or more missing terms by any constant and then apply the boundary conditions which eliminates the constants. The crisp solution of fuzzy boundary value problem (FBVP) always lies between the upper and lower solutions. 
 If the lower solution is not monotonically increasing and the upper solution is not monotonically decreasing then the solution of the FDE is not a valid level set.

\noindent This paper is organized as follows:\\
In section $2$, we recall some basics definitions and theorems. FLT is defined in section $3$ and in this section the FBVP is briefly reviewed. In section $4$, constructing solution of FBVP by FLT is explained. To illustrate the  method, several examples are given in section $5$. Conclusion is given in section $6$.

\section{Basic concepts}
In this section we will recall some basics definitions and theorems needed throughout the paper such as fuzzy number, fuzzy-valued function and the derivative of the fuzzy-valued functions.
\begin{definition} A fuzzy number is defined in \cite{1} as the mapping such that $u:R\rightarrow[0,1]$, which satisfies the following four properties
\begin{enumerate}
\item $u$ is upper semi-continuous.
 \item $u$ is fuzzy  convex that is $u(\lambda  x+(1-\lambda)y) \geq \min{\{u(x), u(y)\}}, x, y\in R$ and $\lambda\in [0,1]$.
\item $u$ is normal that is $\exists$ $x_0\in R$, where $u(x_0)=1$.
\item $A=\{\overline{x \in \mathbb{R}: u(x)>0}\}$ is compact, where $\overline{A}$ is closure of $A$.
\end{enumerate}
\end{definition}

\begin{definition}
A fuzzy number in parametric form given in \cite{2,3,4} is an order pair of the form $u=(\underline{u}(r), \overline{u}(r))$, where $0\leq r\leq1$ satisfying the following conditions.
\begin{enumerate}
\item $\underline{u}(r)$ is a bounded left continuous increasing function in the interval $[0,1]$.
\item $\overline{u}(r)$ is a bounded left continuous decreasing function in the interval $[0,1]$.
\item $\underline{u}{(r)\leq\overline{u}(r)}$. 
\end{enumerate}
If $\underline{u}(r)=\overline{u}(r)=r$, then $r$ is called crisp number.
\end{definition}

Now we recall a triangular fuzzy from \cite{1,18,19} number which must be in the form of
$u=(l, c, r),$ where $l,c,r\in R$ and $l\leq c\leq r$, then $\underline{u}(\alpha)=l+(c-r)\alpha$ and $\overline{u}(\alpha)=r-(r-c)\alpha$ are the end points of the $\alpha$ level set.
Since each $y\in R$ can be regarded as a fuzzy number if
\begin{eqnarray*}\widetilde{y}(t)=\begin{cases}1, \;\;\; if \;\; y=t,\\  0, \;\;\; if \;\; 1\neq t.\end{cases}\end{eqnarray*}
For arbitrary fuzzy numbers $u=(\underline{u}(\alpha), \overline{u}(\alpha))$ and $v=(\underline{v}(\alpha), \overline{v}(\alpha))$ and an arbitrary crisp number $j$, we define addition and scalar multiplication as:
\begin{enumerate}
\item $(\underline{u+v})(\alpha)=(\underline{u}(\alpha)+\underline{v}(\alpha))$.
\item $(\overline{u+v})(\alpha)=(\overline{u}(\alpha)+\overline{v}(\alpha))$.
\item $(j\underline{u})(\alpha)=j\underline{u}(\alpha)$, $(j\overline{u})(\alpha)=j\overline{u}(\alpha)$ \mbox{       }  $j\geq0$.
\item $(j\underline{u})(\alpha)=j\overline{u}(\alpha)\alpha, (j\overline{u})(\alpha)=j\underline{u}(\alpha)\alpha$, $j<0$.
\end{enumerate}
\begin{definition} (See Salahshour \& Allahviranloo, and Allahviranloo \& Barkhordari \cite{18,19}) Let us suppose that x, y $\in E$, if $\exists$ $z\in E$ such that
$x=y+z$. Then, $z$ is called the H-difference of $x$ and $y$ and is given by $x\ominus y$.\end{definition}
\begin{remarks}(see Salahshour \& Allahviranloo \cite{18}).
Let $X$ be a cartesian product of the universes, $X_1$, $X_1, \cdots, X_n$, that is
$X=X_1 \times X_2 \times \cdots \times X_n$ and $A_{1},\cdots,A_{n}$ be $n$ fuzzy numbers in $X_1, \cdots, X_n$ respectively. Then,  $f$ is a mapping from $X$ to a universe $Y$, and $y=f(x_{1},x_{2},\cdots,x_{n})$, then the Zadeh extension principle allows us to define a fuzzy set $B$ in $Y$ as;
\begin{equation*}B=\{(y, u_B(y))|y=f(x_1,\cdots,x_{n}),(x_{1},\cdots,x_{n})\in X\},\end{equation*}
\noindent where
\begin{eqnarray*}u_B(y)=\begin{cases}
\sup_{{(x_1,\cdots,x_n)} \in f^{-1}(y)}  \min\{u_{A_1}(x_1),\cdots u_{A_n}(x_n)\}, \;\;\;  if  \;\;\; f^{-1}(y)\neq 0,\\ 0, \;\;\;\; otherwise,
\end{cases}\end{eqnarray*}
\noindent where $f^{-1}$ is the inverse of $f$.

The extension principle reduces in the case if $n=1$ and is given as follows: 
$B=\{(y, u_B(y)|y=f(x), \mbox{    } x \in X)\},$
\noindent where
\begin{eqnarray*}u_B(y)=\begin{cases}\sup_{x\in f^{-1}(y)} \{u_A(x)\}, \mbox{   if   } f^{-1}(y)\neq 0,\\0, \;\;\;\; otherwise. \end{cases} \end{eqnarray*}

By Zadeh extension principle the approximation of addition of $E$ is defined by
$(u\oplus v)(x)=\sup_{y\in R}  \min(u(y), v(x-y))$, $x \in R$  and scalar multiplication of a fuzzy number is defined by
\begin{eqnarray*}(k\odot u)(x)=\begin{cases}u(\frac{x}{k}), \;\;\; k > 0,\\ 0 \;\;\; \mbox{ otherwise },  \end{cases} \end{eqnarray*}
\noindent where $\widetilde{0}\in E$.
\end{remarks}
\noindent The Housdorff distance between the fuzzy numbers \cite{6,12,18,19} defined by
\[d:E\times E\longrightarrow R^{+}\cup \{{0}\},\]
\[d(u,v)=\sup_{r\in[0,1]}\max\{|\underline{u}(r)-\underline{v}(r)|, |\overline{u}(r)-\overline{v}(r)|\},\] \noindent where $u=(\underline{u}(r), \overline{u}(r))$ and $v=(\underline{v}(r), \overline{v}(r))\subset R$.
\\\\
We know that if $d$ is a metric in $E$, then it will satisfies the following properties, introduced by Puri and Ralescu \cite{28}:
\begin{enumerate}
\item $d(u+w,v+w)=d(u,v)$, $\forall$  u, v, w $\in$ E.

\item $(k \odot u, k \odot v)=|k|d(u, v)$, $\forall$ k $\in$ R, \mbox{  and  } u, v $\in$ E.

\item $d(u \oplus v, w \oplus e)\leq d(u,w)+d(v,e)$, $\forall$ u, v, w, e $\in$  E.
\end{enumerate}
\begin{definition}(see Song and Wu \cite{29}).
If $f:R\times E \longrightarrow E$, then $f$ is continuous at point $(t_0,x_0) \in R \times E$ provided that for any fixed number $r \in [0,1]$ and any $\epsilon > 0$, $\exists$ $\delta(\epsilon,r)$ such that
$d([f(t,x)]^{r}, [f(t_{0},x_{0})]^{r}) < \epsilon$
whenever $|t-t_{0}|<\delta (\epsilon, r)$ and $d([x]^{r}, [x_{0}]^{r})<\delta(\epsilon,r)$ $\forall$ t $\in$ R, x $\in E$.
\end{definition}
\begin{theorem} (see Wu \cite{30}).
Let $f$ be a fuzzy-valued function on $[a,\infty)$ given in the parametric form as $(\underline{f}(x,r), \overline{f}(x,r))$ for any constant number $r\in[0,1]$. Here we assume that $\underline{f}(x,r)$ and $\overline{f}(x,r)$ are Riemann-Integral on $[a,b]$ for every $b\geq a$. Also we assume that $\underline{M}(r)$ and $\overline{M}(r)$ are two positive functions, such that\\
$\int_a^b|\underline{f}(x,r)| dx \leq \underline{M}(r)$ and $\int_a^b |\overline{f}(x,r)| dx \leq \overline{M}(r)$
for every $b\geq a$, then $f(x)$ is improper integral on $[{a}, \infty)$. Thus an improper integral will always be a fuzzy number.\\
In short \[ \int_a^r f(x) dx = ( \int_a^b|\underline{f}(x,r)| dx, \int_a^b |\overline{f}(x,r)| dx).\]
It is will known that Hukuhare differentiability for fuzzy function was introduced by Puri \& Ralescu in $1983$. 
\end{theorem}
\begin{definition}(see Chalco-Cano and Rom\'{a}n-Flores \cite{31}).
Let $f:(a,b)\rightarrow E$ where $x_{0}\in (a,b)$. Then, we say that $f$ is strongly generalized differentiable at $x_0$ (Beds and Gal differentiability).
If $\exists$ an element $f'(x_0)\in E$ such that
\begin{enumerate}
 \item $\forall h>0$ sufficiently small $\exists$ $f(x_0+h)\ominus f(x_0)$, $f(x_0)\ominus f(x_0-h)$, then the following limits hold (in the metric $d$)\\
 $\lim_{h\rightarrow 0}\frac{f(x_0+h)\ominus f(x_0)}{h}=\lim_{h\rightarrow 0}\frac{f(x_0)\ominus f(x_0-h)}{h}=f'(x_0)$,

\noindent Or
\item $\forall h>0$ sufficiently small, $\exists$ $f(x_0)\ominus f(x_0+h)$,  $f(x_0-h)\ominus f(x_0)$, then the following limits hold (in the metric $d$) \\$\lim_{h\rightarrow 0}\frac{f(x_0)\ominus f(x_0+h)}{-h}=\lim_{h\rightarrow 0}\frac{f(x_0-h)\ominus f(x_0)}{-h}=f'(x_0)$,

    \noindent Or
 \item $\forall h>0$ sufficiently small $\exists$ $f(x_0+h)\ominus f(x_0)$,  $f(x_0-h)\ominus f(x_0)$ and the following limits hold (in metric $d$)\\
$\lim_{h\rightarrow 0}\frac{(x_0+h)\ominus f(x_0)}{h}=\lim_{h\rightarrow 0}\frac{f(x_0-h)\ominus f(x_0)}{-h}=f'(x_0)$,

\noindent Or
\item $\forall h>0$ sufficiently small $\exists$ $f(x_0)\ominus f(x_0+h)$,  $f(x_0)\ominus f(x_0-h)$, then the following limits holds(in metric $d$)\\
$\lim_{h\rightarrow 0}\frac{f(x_0)\ominus f(x_0+h)}{-h}=\lim_{h\rightarrow 0}\frac{f(x_0-h)\ominus f(x_0)}{h}=f'(x_0)$.
\end{enumerate}
The denominators $h$ and $-h$ denote multiplication by $\frac{1}{h}$ and $\frac{-1}{h}$ respectively.\end{definition}
\begin{theorem}(Ses Chalco-Cano and Rom\'{a}n-Flores \cite{31}).

Let $f:R\rightarrow E$ be a function denoted by $f(t)=(\underline{f}(t,r),\overline{f}(t,r))$ for each $r\in[0,1]$. Then
\begin{enumerate}
\item If $f$ is $(i)$-differentiable, then $\underline{f}(t,r)$ and $\overline{f}(t,r)$ are differentiable functions and $f'(t)=(\underline{f}'(t,r), \overline{f}'(t,r))$,
\item If $f$ is $(ii)$-differentiable, then $\underline{f}(t,r)$ and $\overline{f}(t,r)$ are differentiable functions and $f'(t)=(\overline{f}'(t,r), \underline{f}'(t,r))$.
\end{enumerate}
\end{theorem}
\begin{lemma}(see Bede and Gal \cite{32,33}).
Let $x_0\in R$. Then, the FDE $y'=f(x,y)$, $y(x_0)=y_0\in R$ and $f:R\times E\rightarrow E$ is supposed to be a continuous and equivalent to one of the following integral equations.
\[y(x)=y_0+\int_{x_0}^x f(t, y(t))dt \;\;\; \forall  \;\;\; x\in [x_0, x_1],\]
\noindent or
\[y(0)=y^1(x)+(-1)\odot\int_{x_0}^x f(t,y(t))dt\;\;\; \forall \;\;\; x\in [x_0, x_1],\]
\noindent on some interval $(x_0, x_1)\subset R$ depending on the strongly generalized differentiability. Integral equivalency shows that if one solution satisfies the given equation, then the other will also satisfy.
\end{lemma}
\begin{remarks}(see Bede and Gal \cite{32,33}).
In the case of strongly generalized differentiability to the FDE's $y'=f(x,y)$ we use two different integral equations. But in the case of differentiability as the definition of H-derivative, we use only one integral. The second integral equation as in Lemma $2.10$ will be in the form of $y^{1}(t)=y^{1}_0\ominus(-1)\int_{x_0}^x f(t,y(t))dt$. The following theorem related to the existence of solution of FIVP under the generalized differentiability.
\end{remarks}

\begin{theorem}
Let us suppose that the following conditions are satisfied.
\begin{enumerate}
\item Let $R_0=[x_0, x_0+s]\times B(y_0, q), s,q>0, y\in E$, where $B(y_0,q)=\{y\in E: B(y,y_0)\leq q\}$ which denotes a closed ball in $E$ and let $f:R_0\rightarrow E$ be continuous functions such that $D(0, f(x,y))\leq M$, $\forall (x,y) \in R_0$ and $0\in E$.
\item Let $g:[x_0, x_0+s]\times [0,q]\rightarrow R$ such that $g(x, 0)\equiv 0$ and $0\leq g(x,u)\leq M$, $\forall x \in [x_0, x_0+s], 0\leq u\leq q$, such that $g(x,u)$ is increasing in u, and g is such that the FIVP $u'(x)=g(x, u(x)), u(x)\equiv 0$ on $[x_0, x_0+s].$
\item  We have $D[f(x,y),f(x,z)\leq g(x, D(y,z))]$,$\forall$ (x,y), (x, z)$\in R_0$ and $D(y,z)\leq q.$

\item $\exists d>0$ such that for $x\in [x_0, x_0+d]$, the sequence $y^1_n:[x_0, x_0+d]\rightarrow E$ given by $y^1_0(x)=y_0$, $y^1_{n+1}(x)=y_0\ominus(-1)\int_{x_0}^{x} f(t, y^{1}_n)dt$ defined for any $n\in N$. Then the FIVP $y'=f(x,y)$, $y(x_0)=y_0$ has two solutions that is (1)-differentiable and the other one is (2)-differentiable for $y$.\end{enumerate}
$y^{1}=[x_0, x_0+r]\rightarrow B(y_0, q)$, where $r=\min\{s,\frac{q}{M},\frac{q}{M_1},d\}$ and the successive iterations $y_0(x)=y_0$, $y_{n+1}(x)=y_{0}+\int_{x_0}^{x}f(t,y_{n}(t))dt$ and $y^{1}_{n+1}=y_0$, $y^{1}_{n+1}(x)=y_0\ominus (-1)\int_{x_0}^{x}f(t, y^{1}_{n}(t))dt$ converge to these two solutions respectively. Now according to theorem (2.11), we restrict our attention to function which are (1) or (2)-differentiable on their domain except on a finite number of points as discussed in \cite{33}.
\end{theorem}
\section{Fuzzy Laplace Transform}
Suppose that $f$ is a fuzzy-valued function and $p$ is a real parameter, then according to \cite{18,19} FLT of the function $f$ is defined as follows:
\begin{definition}\label{eq1}
The FLT of fuzzy-valued function is \cite{18,19}
\begin{equation}\label{eq2}\widehat{F}(p)=L[f(t)]=\int_{0}^{\infty}e^{-pt}f(t)dt,\end{equation}
\begin{equation}\label{eq3}\widehat{F}(p)=L[f(t)]=\lim_{\tau\rightarrow\infty}\int_{0}^{\tau}e^{-pt}f(t)dt,\end{equation}
\begin{equation}\widehat{F}(p)=[\lim_{\tau\rightarrow\infty}\int_{0}^{\tau}e^{-pt}\underline{f}(t)dt,\lim_{\tau\rightarrow\infty}\int_{0}^{\tau}e^{-pt}\overline{f}(t)dt],\end{equation}
\noindent whenever the limits exist.
\end{definition}
\begin{definition}\textbf{Classical Fuzzy Laplace Transform:} Now consider the fuzzy-valued function in which the lower and upper FLT of the function are represented by
\begin{equation}\label{eq4}\widehat{F}(p;r)=L[f(t;r)]=[l(\underline{f}(t;r)),l(\overline{f}(t;r))],\end{equation}
\noindent where
\begin{equation}\label{eq5}l[\underline{f}(t;r)]=\int_{0}^{\infty}e^{-pt}\underline{f}(t;r)dt=\lim_{\tau\rightarrow\infty} \int_{0}^{\tau}e^{-pt}\underline{f}(t;r)dt,\end{equation}
\begin{equation}\label{eq6}l[\overline{f}(t;r)]=\int_{0}^{\infty}e^{-pt}\overline{f}(t;r)dt=\lim_{\tau\rightarrow\infty}\int_{0}^{\tau}e^{-pt}\overline{f}(t;r)dt.  \end{equation}
\end{definition}
\subsection{Fuzzy Boundary Value problem}
The concept of fuzzy numbers and fuzzy set was first introduced by Zadeh \cite{1}. Detail information of fuzzy numbers and fuzzy arithmetic can be found in \cite{12,13,14}. In this section we review the fuzzy boundary valued problem (FBVP) with crisp linear differential equation but having fuzzy boundary values. For example we consider the second order fuzzy boundary problem as \cite{9,10,21,1a}.

\begin{equation}\begin{split}\label{eq76}
\psi^{''}(t)+c_1(t)\psi^{'}(t)+c_2(t)\psi(t)=f(t),\\
\psi(0)=\tilde{A},\\
\psi(l)=\tilde{B}.
\end{split}\end{equation}

\section{Constructing Solutions Via FBVP}
In this section we consider the following second order FBVP in general form under generalized H-differentiability proposed in \cite{21}. We define
\begin{equation}\label{eq7711a}y''(t)=f(t, y(t),y'(t)), \end{equation}
\noindent subject to two-point boundary conditions \\ \[y(0)=(\underline{y}(0;r), \overline{y}(0;r)),\] \[y(l)=(\underline{{y}}(l;r), \overline{y}(l;r)).\]

\noindent Taking FLT of (\ref{eq7711a})
\begin{equation}\label{eq78}L[y''(t)]=L[f(t, y(t),y'(t))], \end{equation}
\noindent which can be written as
\[p^2L[y(t)]\ominus py(0)\ominus y'(0)=L[f(t, y(t),y'(t))].\]
The classical form of FLT is given below:
\begin{equation}\begin{split}\label{eq79}
p^{2}l[\underline{y}(t;r)]-p\underline{y}(0;r)-\underline{y}'(0;r)=l[\underline{f}(t, y(0;r),y'(0;r))],
\end{split}\end{equation}
\begin{equation}\begin{split}\label{eq80}
p^{2}l[\overline{y}(t;r)]-p\overline{y}(0;r)-\overline{y}'(0;r)=l[\overline{f}(t, y(0;r),y'(0;r))].
\end{split}\end{equation}
Here we have to replace the unknown value $y'(0,r)$ by  constant $F_1$  in lower case and by  $F_2$ in upper case. Then we can find these values by applying the given boundary conditions.

In order to solve equations (\ref{eq79}) and (\ref{eq80}) we assume that $A(p;r)$ and $B(p;r)$ are the solutions of (\ref{eq79}) and (\ref{eq80}) respectively. Then the above system becomes
\begin{equation}\label{eq81}
l[\underline{y}(t;r)]=A(p;r),
\end{equation}
\begin{equation}\label{eq30}
l[\overline{y}(t;r)]=B(p;r).
\end{equation}
\noindent Using inverse Laplace transform, we get the upper and lower solution for given problem as:
\begin{equation}\label{eq31}
[\underline{y}(t;r)]=l^{-1}[A(p;r)],
\end{equation}
\begin{equation}\label{eq32}
[\overline{y}(t;r)]=l^{-1}[B(p;r)].
\end{equation}
\section{Examples}
In this section first we consider the $Schr\ddot{o}dinger$ equation \cite{21} with fuzzy boundary conditions under Hukuhara differentiability.
\begin{example}
The $Schr\ddot{o}dinger$ FBVP \cite{21} is as follows:
\begin{equation}\label{eq1n}
(\frac{h^2}{2m})u^{''}(x)+V(x)u(x)=Eu(x),
\end{equation}
\noindent where $V(x)$ is potential and is defined as \begin{eqnarray*}V(x)=\begin{cases}0, \;\;\; if \;\; x < 0,\\  l, \;\;\; if \;\; x>0,\end{cases}\end{eqnarray*}
\noindent subject to the following boundary conditions\\
\[u(0)=(1+r, 3-r),\] \[u(l)=(4+r, 6-r).\]\\ \noindent Now let $a=\frac{h^2}{2m}$, $b=E$. Then, (\ref{eq1n}) becomes
\begin{equation}\label{eq2n}
au^{''}(x)+V(x)u(x)=bu(x).
\end{equation}
In (\ref{eq1n}) for $x<0$, we discuss (1,1) and (2,2)-differentiability while in the case $x>0$ we will discuss (1,2) and (2,1)-differentiability.
\subsection{Case-I: (1,1) and (2,2)-differentiability}
For $x<0$, (\ref{eq2n}) becomes
\begin{equation*}
au^{''}=bu,
\end{equation*}
\begin{equation}\label{eq4n}
au^{''}-bu=0.
\end{equation}
Now applying FLT on both sides of equation (\ref{eq4n}), we get
\begin{equation*}\label{eq65}
aL[u^{''}(x)]-bL[u(x)]=0,
\end{equation*}
\noindent where
\begin{equation*}\label{eq73}
L[u''(x)]=p^{2}L[u(x)]\ominus pu(0)\ominus u'(0).
\end{equation*}
The classical FLT form of the above equation is
\begin{equation*}\label{eq74}
l[\underline u''(x,r)]=p^{2}l[\underline{u}(x,r)]-p\underline u(0,r)-\underline u'(0,r),
\end{equation*}
\begin{equation*}\label{eq75}
l[\overline u''(x,r)]=p^{2}l[\overline {u}(x,r)]-p\overline u(0,r)-\overline u'(0,r).
\end{equation*}
\noindent Solving  the above classical equations for lower and upper solutions, we have
\begin{equation*}\begin{split}\label{eq76}
a\{p^{2}l[\underline{u}(x,r)]-p\underline u(0,r)-\underline u'(0,r)\}-bl[\underline {u}(x,r)]=0,
\end{split}\end{equation*}
or \begin{equation*}\begin{split}\label{eq76}
(ap^{2}-b)l[\underline{u}(x,r)]=a\{p\underline u(0,r)+\underline u'(0,r)\}.
\end{split}\end{equation*}
\noindent Applying the boundary conditions, we have
\begin{equation*}\begin{split}\label{eq77}
(ap^{2}-b)l[\underline{u}(x,r)]=a\{p(1+r)+F_1\},
\end{split}\end{equation*}
where \[F_1=\underline u'(0,r).\]
\noindent Simplifying and applying inverse Laplace we get
\begin{equation*}\begin{split}\label{eq77}
\underline{u}(x,r)=(\frac{1+r}{2})l^{-1}\{\frac{p}{p^2-\frac{b}{a}}\}+F_1l^{-1}\{\frac{1}{p^2-\frac{b}{a}}\}.
\end{split}\end{equation*}
Using partial fraction
\begin{equation}\begin{split}\label{eq7}
\underline{u}(x,r)=(\frac{1+r}{2})\{e^{\sqrt{\frac{b}{a}}x}+e^{-\sqrt{\frac{b}{a}}x}\}+\frac{F_1}{2\sqrt{\frac{b}{a}}}\{e^{\sqrt{\frac{b}{a}}x}-e^{-\sqrt{\frac{b}{a}}x}\}.
\end{split}\end{equation}
\noindent Now applying boundary conditions on (\ref{eq7}) we get
\begin{equation*}\begin{split}\label{eq77}
F_1=\frac{4+r-\frac{1+r}{2}\{e^{\sqrt{\frac{b}{a}}l}+e^{-\sqrt{\frac{b}{a}}l}\}}{\frac{1}{2}\sqrt{\frac{a}{b}}\{e^{\sqrt{\frac{b}{a}}l}-e^{-\sqrt{\frac{b}{a}}l}\}}.
\end{split}\end{equation*}
\noindent Putting value of $F_1$ in (\ref{eq7}) we get
\begin{equation*}\begin{split}\label{eq771}
\underline{u}(x,r)=(\frac{1+r}{2})\{e^{\sqrt{\frac{b}{a}}x}+e^{-\sqrt{\frac{b}{a}}x}\}+\frac{4+r-\frac{1+r}{2}\{e^{\sqrt{\frac{b}{a}}l}+e^{-\sqrt{\frac{b}{a}}l}\}}{\{e^{\sqrt{\frac{b}{a}}l}-e^{-\sqrt{\frac{b}{a}}l}\}}\{e^{\sqrt{\frac{b}{a}}x}-e^{-\sqrt{\frac{b}{a}}x}\}.
\end{split}\end{equation*}
Now solving the classical FLT form for $\overline{u}(x,r)$, we have
 \begin{equation*}\begin{split}\label{eq76}
a\{p^{2}l[\overline {u}(x,r)]-p\overline  u(0,r)-\overline u'(0,r)\}-bl[\overline  {u}(x,r)]=0,
\end{split}\end{equation*}
 \begin{equation*}\begin{split}\label{eq76}
(ap^{2}-b)l[\overline {u}(x,r)]=a\{p\overline  u(0,r)+\overline  u'(0,r)\}.
\end{split}\end{equation*}
\noindent Using the boundary conditions, we have
\begin{equation*}\begin{split}\label{eq77}
(ap^{2}-b)l[\overline {u}(x,r)]=a\{p(3-r)+F_2\},
\end{split}\end{equation*}
\noindent where \[F_2=\overline  u'(0,r).\]
\noindent Simplifying and applying inverse laplace we get
\begin{equation*}\begin{split}\label{eq77}
\overline {u}(x,r)=(\frac{3-r}{2})l^{-1}\{\frac{p}{p^2-\frac{b}{a}}\}+F_2l^{-1}\{\frac{1}{ap^2-\frac{b}{a}}\}.
\end{split}\end{equation*}
Using partial fraction
\begin{equation}\begin{split}\label{eq774}
\overline{u}(x,r)=(\frac{3-r}{2})\{e^{\sqrt{\frac{b}{a}}x}+e^{-\sqrt{\frac{b}{a}}x}\}+\frac{F_2}{2\sqrt{\frac{b}{a}}}\{e^{\sqrt{\frac{b}{a}}x}-e^{-\sqrt{\frac{b}{a}}x}\}.
\end{split}\end{equation}
\noindent Now applying boundary conditions on (\ref{eq774}) we have
\begin{equation*}\begin{split}\label{eq77}
F_2=\frac{6-r-\frac{3-r}{2}\{e^{\sqrt{\frac{b}{a}}l}+e^{-\sqrt{\frac{b}{a}}l}\}}{\frac{1}{2}\sqrt{\frac{a}{b}}\{e^{\sqrt{\frac{b}{a}}l}-e^{-\sqrt{\frac{b}{a}}l}\}}.
\end{split}\end{equation*}
\noindent Putting value of $F_2$ in (\ref{eq774}) we get
\begin{equation*}\begin{split}\label{eq771}
\overline{u}(x,r)=(\frac{3-r}{2})\{e^{\sqrt{\frac{b}{a}}x}+e^{-\sqrt{\frac{b}{a}}x}\}+\frac{6-r-\frac{3-r}{2}\{e^{\sqrt{\frac{b}{a}}l}+e^{-\sqrt{\frac{b}{a}}l}\}}{\{e^{\sqrt{\frac{b}{a}}l}-e^{-\sqrt{\frac{b}{a}}l}\}}\{e^{\sqrt{\frac{b}{a}}x}-e^{-\sqrt{\frac{b}{a}}x}\}.
\end{split}\end{equation*}
\subsection{Case-II: (1) and (2)-differentiability, (2) and (1)-differentiability}
For $x>0$, (\ref{eq2n}) becomes
\begin{equation}\label{eq4nn}
au^{''}+(l-b)u=0.
\end{equation}
Applying FLT and inverse Laplace transform and then simplifying  we get the following lower solution.
 \begin{equation*}\begin{split}\label{eq771}
\underline{u}(x,r)=(\frac{1+r}{2})\bigg[\cos\frac{x\sqrt{b-l}}{\sqrt{a}}+\cosh\frac{x\sqrt{b-l}}{\sqrt{a}}\bigg]+\frac{H_1\sqrt{a}}{2\sqrt{b-l}}\bigg[\sin\frac{x\sqrt{b-l}}{\sqrt{a}}+\sinh\frac{x\sqrt{b-l}}{\sqrt{a}}\bigg]\\
-\frac{(3-r)}{2}\bigg[\cos\frac{x\sqrt{b-l}}{\sqrt{a}}-\cosh\frac{x\sqrt{b-l}}{\sqrt{a}}\bigg]-\frac{H_2\sqrt{a}}{2\sqrt{b-l}}\bigg[\sin\frac{x\sqrt{b-l}}{\sqrt{a}}-\sinh\frac{x\sqrt{b-l}}{\sqrt{a}}\bigg],
\end{split}\end{equation*}
\noindent or
\begin{equation*}\begin{split}\label{eq771}
\underline{u}(x,r)=\frac{1+r}{2}(c_1)+\frac{H_1\sqrt{a}}{2\sqrt{b-l}}(c_2)
-\frac{3-r}{2}(c_3)-\frac{H_2\sqrt{a}}{2\sqrt{b-l}}(c_4).
\end{split}\end{equation*}

\noindent The upper solution will be as follows:
\begin{equation*}\begin{split}\label{eq771}
\overline{u}(x,r)=(\frac{3-r}{2})\bigg[\cos\frac{x\sqrt{b-l}}{\sqrt{a}}+\cosh\frac{x\sqrt{b-l}}{\sqrt{a}}\bigg]+\frac{H_2\sqrt{a}}{2\sqrt{b-l}}\bigg[\sin\frac{x\sqrt{b-l}}{\sqrt{a}}+\sinh\frac{x\sqrt{b-l}}{\sqrt{a}}\bigg]\\
-\frac{(1+r)}{2}\bigg[\cos\frac{x\sqrt{b-l}}{\sqrt{a}}-\cosh\frac{x\sqrt{b-l}}{\sqrt{a}}\bigg]-\frac{H_1\sqrt{a}}{2\sqrt{b-l}}\bigg[\sin\frac{x\sqrt{b-l}}{\sqrt{a}}-\sinh\frac{x\sqrt{b-l}}{\sqrt{a}}\bigg]
\end{split}\end{equation*}
\begin{equation*}\begin{split}\label{eq771}
\overline{u}(x,r)=\frac{3-r}{2}(c_1)+\frac{H_2\sqrt{a}}{2\sqrt{b-l}}(c_2)
-\frac{1+r}{2}(c_3)-\frac{H_1\sqrt{a}}{2\sqrt{b-l}}(c_4)
\end{split}\end{equation*}
\noindent where \[c_1=\cos\frac{x\sqrt{b-l}}{\sqrt{a}}+\cosh\frac{x\sqrt{b-l}}{\sqrt{a}},\] \[c_2=\sin\frac{x\sqrt{b-l}}{\sqrt{a}}+\sinh\frac{x\sqrt{b-l}}{\sqrt{a}},\] \[c_3=\cos\frac{x\sqrt{b-l}}{\sqrt{a}}-\cosh\frac{x\sqrt{b-l}}{\sqrt{a}},\]  \[c_4=\sin\frac{x\sqrt{b-l}}{\sqrt{a}}-\sinh\frac{x\sqrt{b-l}}{\sqrt{a}}.\]
\[H_1=\frac{2c_2}{c^2_2-c^2_4}\bigg[4+r-\frac{r+1}{2}c_1+\frac{3-r}{2}c_3\bigg]+\frac{2c_4}{c^2_2-c^2_4}\bigg[6-r-\frac{3-r}{2}c_1+\frac{1+r}{2}c_3\bigg],\]
and
\[H_2=\frac{2c_4}{c^2_2-c^2_4}\bigg[4+r-\frac{r+1}{2}c_1+\frac{3-r}{2}c_3\bigg]+\frac{2c_2}{c^2_2-c^2_4}\bigg[6-r-\frac{3-r}{2}c_1+\frac{1+r}{2}c_3\bigg].\]
\end{example}

\begin{example} Consider the following fuzzy homogenous boundary value problem
\begin{equation}\label{s1}
x^{''}(t)-3x^{'}(t)+2x(t)=0,\end{equation}
\noindent subject to the following boundary conditions\\ \[x(0)=(0.5r-0.5, 1-r),\] \[x(1)=(r-1, 1-r).\]\\
Now applying fuzzy Laplace transform on both sides of equation (\ref{s1}), we get
\begin{equation}\label{eq65}
L[x^{''}(t)]=3L[x'(t)]-2L[x(t)].
\end{equation}
\noindent We know that
\begin{equation*}\label{s2}
L[x''(t)]=p^{2}L[x(t)]\ominus px(0)\ominus x'(0).
\end{equation*}
The classical FLT form of the above equation is
\begin{equation*}\label{eq74}
l[\underline x''(t,r)]=p^{2}l[\underline{x}(t,r)]-p\underline{x}(0,r)-\underline{x}'(0,r),
\end{equation*}
\begin{equation*}\label{eq75}
l[\overline x''(t,r)]=p^{2}l[\overline {x}(t,r)]-p\overline{x}(0,r)-\overline x'(0,r).
\end{equation*}
Now on putting in (\ref{eq65}), we have
\begin{equation}\label{eq76a}
p^{2}l[\underline{x}(t,r)]-p\underline x(0,r)-\underline x'(0,r)\}-3pl[\underline {x}(t,r)]+3\underline{x}(0,r)+2l[\underline{x}(t,r)]=0,\end{equation}
\begin{equation}\label{eq76b}p^{2}l[\overline{x}(t,r)]-p\overline{x}(0,r)-\overline x'(0,r)\}-3pl[\overline{x}(t,r)]+3\overline{x}(0,r)+2l[\overline{x}(t,r)]=0.
\end{equation}
Solving (\ref{eq76a}) for $l[\underline{x}(t,r)]$, we get
 \begin{equation*}\begin{split}\label{eq76}
(p^{2}-3p+2)l[\underline{x}(t,r)]=p\underline x(0,r)+\underline x'(0,r)+3[\underline {x}(0,r)].
\end{split}\end{equation*}
\noindent Applying boundary conditions, we have
\begin{equation*}\label{eq76}
l[\underline{x}(t,r)]=\frac{(0.5r-0.5)p}{p^2-3p+2}-\frac{3(0.5r-0.5)}{p^2-3p+2}+\frac{A}{p^2-3p+2}.
\end{equation*}
Using partial fraction and then applying inverse Laplace, we get
\begin{equation*}
 \underline{x}(t,r)=(0.5r-0.5)[-e^t+2e^{2t}]-3(0.5r-0.5)[-e^t+2e^{2t}]+A[-e^t+e^{2t}].\end{equation*}
 Using boundary values, we get
 \begin{equation*}\underline{x}(1,r)=r-1=(0.5r-0.5)[-e+2e^2]-3(0.5r-0.5)[-e+2e^2]+A[-e+e^2],\end{equation*}

 \begin{equation*}A=\frac{r-1+(0.5r-0.5)[-2e+e^2]}{e^2-e}.\end{equation*}
 \noindent Finally on putting value of A we have
 \begin{equation*}
  \underline{x}(t,r)=(0.5r-0.5)(-e^t+2e^{2t})-3(0.5r-0.5)(-e^t+2e^{2t})+\frac{r-1+(0.5r-0.5)(-2e+e^2)}{e^2-e}(-e+e^2) \end{equation*}
Now solving (\ref{eq76b}) for $l[\overline{x}(t,r)]$, we have
\begin{equation*}\label{eq76}
(p^{2}-3p+2)l[\overline{x}(t,r)]=p\overline x(0,r)+\overline x'(0,r)\}+3[\overline {x}(0,r)].
\end{equation*}
\noindent Applying boundary condition we get
\begin{equation*}\label{eq76}
l[\overline{x}(t,r)]=\frac{(1-r)p}{p^2-3p+2}-\frac{3(1-r)}{p^2-3p+2}+\frac{A}{p^2-3p+2}.
\end{equation*}
\noindent Using partial fraction and then applying inverse Laplace
\begin{equation}\label{eq76c}
 \overline{x}(t,r)=(1-r)[-e^t+2e^{2t}]-3(1-r)[-e^t+2e^{2t}]+A[-e^t+e^{2t}].\end{equation}
 Using boundary values
 \begin{equation*}
 \overline{x}(1,r)=1-r=(1-r)[-e+2e^2]-3(1-r)[-e+2e^2]+A[-e+e^2],\end{equation*}
 \begin{equation*}A=\frac{1-r+(1-r)[-2e+e^2]}{e^2-e}.\end{equation*}
\noindent Putting value of A in (\ref{eq76c}) we get
\begin{equation*}
\overline {x}(t,r)=(1-r)[-e^t+2e^{2t}]-3(1-r)[-e^t+2e^{2t}]+\frac{1-r+(1-r)[-2e+e^2]}{e^2-e}[-e+e^{2t}].
\end{equation*}
\end{example}

\section{Conclusion}
In this paper, we applied the fuzzy Laplace transform to solve FBVPs under generalized H-differentiability, in particular,
solving $Schr\ddot{o}dinger$ FBVP. We also used FLT to solve homogenous FBVP. This is another application of FLT. Thus FLT can also be used to solve FBVPs analytically. The method can be extended for an $nth$ order FBVP. This work is in progress.

\end{document}